\documentclass[10pt]{article}

\usepackage{amssymb}
\usepackage[leqno]{amsmath}
\usepackage{amsthm}

\newtheorem{theorem}{Theorem}[section]
\newtheorem{corollary}[theorem]{Corollary}
\newtheorem{proposition}[theorem]{Proposition}
\newtheorem{lemma}[theorem]{Lemma}

\theoremstyle{definition}

\newtheorem{remark}[theorem]{Remark}

\newcommand\Symp{\operatorname{Symp}}

\newcommand\Def{{\overset {\rm {def}}{\ =\ }}}

\def\R{\mathbb R}
\def\C{\mathbb C}

\def\Z{\mathbb Z}

\def\F{\mathcal F}

\def\J{\mathcal J}
\def\M{\mathcal M}

\def\S{\mathcal S}
\def\B{\mathcal B}

\begin{document}

\bibliographystyle{abbrv}

\title{Ruled $4$-manifolds and \\ isotopies of symplectic surfaces}

\author{R. Hind \thanks{Supported in part by NSF grant DMS-0505778.}
\and A. Ivrii
}

\date{\today}
\maketitle

\abstract We study symplectic surfaces in ruled symplectic
$4$-manifolds which are disjoint from a given symplectic section.
As a consequence, in any symplectic $4$-manifold, two homologous
symplectic surfaces which are $C^0$ close must be Hamiltonian
isotopic.

\endabstract

\section{Introduction}

This paper's initial objective was to address a question of local isotopy
of symplectic surfaces in symplectic $4$-manifolds (and this was originally
motivated by the question of local isotopy of lagrangian surfaces in $4$-manifolds
\cite{HI}). More precisely, let
$(M, \omega)$ be a symplectic $4$-manifold and let $\Sigma\subset M$ be a compact embedded
symplectic surface. It is natural to ask whether other symplectic surfaces in some
neighborhood of $\Sigma$ are smoothly or symplectically isotopic to $\Sigma$
(that is isotopic through a family of smooth or symplectic embedded surfaces).
In fact we will show (see corollary~\ref{main_theorem_local}) the following:
\begin{theorem}\label{intro_local}
The set of symplectic
surfaces homologous to $\Sigma$ and lying in a standard symplectic neighborhood
of $\Sigma$ is (weakly) contractible.
\end{theorem}
\noindent
In section $2$ we will briefly describe what these
standard symplectic neighborhoods (standard symplectic disk bundles over $\Sigma$) are,
and refer the reader to \cite{Biran} for a more
complete description of this.
We point out that any symplectic isotopy $\Sigma_t$, $t\in [0,1]$,
is given by a Hamiltonian whose support at each $t$ is arbitrarily close to
$\Sigma_t$ (a proof of this can be found in, for example, \cite{Seibert_Tian}).
We would like to add that in general homologous symplectic surfaces
inside a symplectic $4$-manifold do not need to be
isotopic (see \cite{Fintushel_Stern}), however there are also cases when
this is true (\cite{Shevchishin, Seibert_Tian}).

This paper's real objective is to investigate the set of symplectic surfaces in ruled
symplectic $4$-manifolds. A $4$-manifold $\bar X$ is called ruled if it is an $S^2$-fibration
over a Riemann surface $\Sigma$, such symplectic manifolds have been extensively
studied in the literature (see, for example, Lalonde-McDuff \cite{Lalonde_McDuff1}, \cite{Lalonde_McDuff2},
and also \cite{Abreu_McDuff}, \cite{McDuff01}).
Our main theorem (see theorem~\ref{main_theorem_ruled}) can be formulated as follows:
\begin{theorem}\label{intro_ruled}
Let $Z$ be a symplectic surface in $\bar X$ homologous to a section.
Consider the homology class of sections disjoint from $Z$ and suppose that
it has positive symplectic area.
Then the set of symplectic surfaces in $\bar X$ belonging to this homology class
and avoiding $Z$ is non-empty and (weakly) contractible.
\end{theorem}
\noindent
In particular the statement on local isotopy is an easy consequence of this.
We will discuss ruled symplectic $4$-manifolds in more detail in section $3$,
we also postpone till this section the section-by-section description of the rest of the paper.
Similar results have already been obtained by Joseph Coffey (see \cite{Coffey}), and in fact
the current version of this paper is highly motivated by his results.

The group of symplectomorphisms of $\bar X$ acts on the set of symplectic surfaces in $\bar X$
and this leads to an interesting connection between the two. In particular we can obtain
some results in the spirit of \cite{Abreu_McDuff}.
This is briefly discussed in section $9$.

Unless stated otherwise, we will always mean that a surface is embedded, oriented and non-parameterized,
any family is smoothly dependent on parameters, and homotopical properties
such as homotopy equivalence or contractibility are to be understood in the weak
sense (that is, weak homotopy equivalence or weak contractibility).


{\it Acknowledgments}. The second author would like to thank Shengda Hu, Octavian Cornea and Fran{\c c}ois Lalonde
for many helpful discussions on the subject, Joseph Coffey for explaining his relating work and referring us
to his paper,
and Paul Biran and Leonid Polterovich for a chance to present
the previous approach to the problem on Geometry and Dynamics seminar and helpful comments received therein.

\section{Standard symplectic disk and sphere bundles}\label{standard_bundles}

{\bf Standard symplectic disk bundles: }

Given a Riemann surface $\Sigma$ of genus $g$ and $d\in \Z$, we will construct the most
canonical symplectic disk bundle $(X, \omega)$ of degree $d$ over $\Sigma$
(see also \cite{Biran}).

Denote by $L$ the complex line bundle over $\Sigma$ with Chern number $d$,
thus $[\Sigma]\in H_2(L)$ satisfies $[\Sigma]\cdot [\Sigma] =d$.
Let $r$ denote the radial coordinate on $L$ with respect to a Hermitian
metric on $L$,
let $\pi: L\rightarrow \Sigma$ denote the natural projection to the base,
let $\Sigma_0$ denote the $0$-section,
and let $P = \{ r = 1 \}\subset L$ denote the unit circle bundle.

First suppose that $d\ne 0$.
Fix an area form $\sigma$ on $\Sigma$ with $\int_\Sigma \sigma = |d|$.
We can choose an $S^1$-invariant connection $1$-form $\alpha$ on $P$ so that
$d\alpha =  \pi^*\sigma$,
and extend $\alpha$ to the complement of the $0$-section of $L$ as $\pi^*\alpha$.
In the case $d > 0$,
we let $X =\{ r < 1 \}\subset L$ be the unit disk bundle and
define $\omega = \pi^*\sigma - d(r^2 \alpha)$.
In the case $d < 0$, we let $X = L$ and
define $\omega = \pi^*\sigma + d(r^2 \alpha)$.
Note that in both cases on the complement of the $0$-section the symplectic form can be
written as $\omega = d((1\pm r^2) \alpha)$.

In the case $d=0$, we let $X = \Sigma \times \C$ be the trivial line bundle
over $\Sigma$. We fix an area form $\sigma$ on $\Sigma$ of total area $1$,
and let $\tau$ be the standard symplectic form on $\R^2 = \C$.
Define $\omega = \sigma \oplus \tau$.

In all cases, $\omega$ is a symplectic form on $X$, and the fibers of $X$ and
the $0$-section of $X$ are $\omega$-symplectic.

\vskip5pt
\noindent
{\bf Symplectic neighborhood theorem: }

Let $(M, \omega)$ be any $4$-dimensional symplectic manifold and let $\Sigma\subset M$
be a sympectic surface. By a symplectic neighborhood theorem,
a neighborhood of $\Sigma$ is determined by the area
$\int_{\Sigma} \omega$ and the self-intersection number $d = \Sigma \cdot \Sigma$, the latter
being also equal to the first Chern number of the symplectic normal bundle to $\Sigma$.

Thus in a neighborhood $U$ of $\Sigma$ an appropriate multiple $c\omega$
of the symplectic form $\omega$ can be written in a standard form described above:
\begin{equation}\label{standard_form}
c\omega =
\left\{
\begin{array}{lc}
\pi^*\sigma - d(r^2\alpha) = d((1 - r^2)\alpha) & \text{if }d>0 \\
\pi^*\sigma + d(r^2\alpha) = d((1 + r^2)\alpha) & \text{if }d<0\\
\pi^*\sigma + dy_1 \wedge dy_2 & \text{if }d= 0,\\
\end{array}
\right.
\end{equation}
where $\sigma$ is an area form on $\Sigma$ normalized so that
the total area of $\Sigma$ equals $|d|$ if $d\ne 0$ and equals $1$ if $d=0$,
$\pi: U \rightarrow \Sigma$ is a projection map, and
$\alpha$ (if $d\ne 0$) is a connection $1$-form
with $d\alpha = \pi^*\sigma$.

Moreover the same theorem with parameters applies to any compact family
of $\omega_\lambda$-symplectic surfaces $\Sigma_\lambda$ to produce
a family of such standard presentations of $\omega_\lambda$ in
sufficiently small neighborhoods of $\Sigma_\lambda$.

\vskip5pt
\noindent
{\bf Compactifying standard disk bundles to standard sphere bundles: }

Given a standard symplectic disk bundle $X$ over $\Sigma$,
let $X(R) = \{ r\le R \}$, where $R < 1$ in the case $d>0$.

If $d\ne 0$ then by symplectic cutting (see for example \cite{Biran} for a proof)
one can compactify
the symplectic disk bundle $(X(R), \omega)$ to a symplectic sphere bundle without changing the
area of the fiber. Slightly abusing notation, we denote this sphere bundle by
$(\bar X(R), \omega)$.
Roughly speaking this construction collapses the boundary of $X(R)$ to a symplectic section $\Sigma_\infty$
of $\bar X$, in particular $(\bar X(R)\setminus \Sigma_\infty, \omega)$ is symplectomorphic to
$(X(R), \omega)$.

If $d=0$, we note that a standard symplectic disk $(D(R), \tau)$ can be compactified to a symplectic sphere,
and so trivial symplectic disk bundle $(X(R), \omega = \sigma\oplus\tau)$ can be compactified
to a trivial symplectic sphere bundle $(\bar X(R), \omega = \sigma\oplus \tau)$.
\noindent
When clear, we omit the parameter $R$ from the notation.

The sphere bundles that we obtain are examples of ruled symplectic $4$-manifolds discussed in the next section.
Note that up to a diffeomorphism there are precisely two sphere bundles over a surface of genus $g$: a trivial
one and a nontrivial one.

\section{Ruled symplectic $4$-manifolds}

We recall that a symplectic $4$-manifold $(\bar X, \omega)$ is called ruled if it is a
fibration over a Riemann surface with fiber $S^2$.
By choosing the right orientations of the base and of the fiber
of the fibration, we will always assume that the integral of
$\omega$ over the fiber and the integral of $\omega^2$ over $\bar X$
are positive.
Up to a symplectomorphism ruled $4$-manifolds have been classified by F. Lalonde and D. McDuff
\cite{Lalonde_McDuff1, Lalonde_McDuff2} who showed that any two cohomologous symplectic forms on $\bar X$
are in fact diffeomorphic.

What makes understanding symplectic properties of ruled $4$-manifolds easier
is that for a large class $\J$ of $\omega$-tame almost complex structures  (all
$\omega$-tame almost complex structures if $g>0$) $\bar X$ admits a
holomorphic foliation by smooth embedded spheres representing the homology class
of the fiber. For a general theory of (pseudo)holomorphic curves we refer the
reader to \cite{Gromov, McDuff_Salamon2}.

We introduce several notations which will be used throughout the paper. We denote by $F$ the
homology class of a fiber and by $A_d$ the homology class of a section with
self-intersection $d\in \Z$.
We note that if up to a diffeomorphism $\bar X$ is a trivial bundle over $\Sigma$
then all homology classes $A_d$ with $d$ even can be represented by embedded surfaces
in $\bar X$, and if $\bar X$ is a nontrivial bundle -- then all homology classes $A_d$
with $d$ odd. We also note that $H_2(\bar X; \R)$ is generated by $F$ and $A_d$
for any $d\in \Z$, and $F\cdot F= 0$, $A_d \cdot A_d = d$, $A_d \cdot F = 1$.
We let $a_d(\omega)$ and $f(\omega)$ denote the integrals of
$\omega$ on the homology classes $A_d$ and $F$ respectively.

For a fixed $d\in \Z$ we will say that $(\bar X, \omega)$ satisfies the {\it cohomology
assumption $(*_d)$} if:
\begin{itemize}
\item either $g\ne 0$ or $d\ne 0$,
\item
$a_d(\omega) > 0$ and $a_{-d}(\omega) > 0$.
\end{itemize}
When the symplectic form is understood from the context, we will just write $a_d$ and $f$.

The standard symplectic sphere bundles of degree $d$ we obtained in the previous section
by compactification automatically
satisfy this assumption if either $d\ne 0$ or $g\ne 0$ since both sections $\Sigma_0$
and $\Sigma_\infty$ are $\omega$-symplectic. For a more general ruled $4$-manifold
this assumption just means that both homology classes $A_d$ and $A_{-d}$ can be represented by
sympectic surfaces.
%

As mentioned in the introduction, the main goal of this paper is to understand the set of
symplectic surfaces in $\bar X$. We denote by $\S_d(\bar X, \omega)$ the set of
$\omega$-symplectic surfaces in $\bar X$ representing the homology class $A_d$,
and we denote by $\S(\bar X, \omega)$ the union of $\S_d(\bar X, \omega)$
over all $d\in \Z$. As follows from the
classification of ruled $4$-manifolds, a necessary and sufficient condition for
$\S_d(\bar X, \omega)$ to be nonempty is that $a_d(\omega)> 0$. When clear from the
context, we will omit the arguments $\bar X$ or $\omega$ from the notation.
Given a symplectic surface $Z\in \S_{-d}(\bar X, \omega)$, we will also
denote by $\S_d(\bar X\setminus Z, \omega)$ the set of symplectic surfaces
in $\S_d(\bar X, \omega)$ disjoint from $Z$.

We describe the contents of later sections of this paper.
In section $4$ we prove a theorem on the existence of holomorphic foliations
of $\bar X$ by spheres in the homology class $F$ -- this is very
well-known if $g\ge 1$, but possibly our version
is slightly uncommon in the case $g=0$.
A direct consequence of this existence is the proof in section $5$
that any $S^n$-family $\Sigma_\lambda \subset \S( \bar X\setminus Z, \omega )$
of symplectic surfaces disjoint from a symplectic surface $Z\in \S_{-d}(\bar X, \omega)$
can be contracted through smoothly embedded surfaces.

Given a symplectic surface $\Sigma\in \S_d$, a standard presentation (\ref{standard_form})
of the symplectic form in a neighborhood of $\Sigma$, and a symplectic foliation $\F$ of $\bar X$
with all leaves intersecting $\Sigma$ positively and transversely, we show in section $6$
how to modify this foliation near $\Sigma$ so that its leaves near $\Sigma$ coincide with the
fibers of $\pi$. We will say that $\F$ {\it nicely intersects} $\Sigma$ to denote this phenomenon.

Many properties of ruled symplectic $4$-manifolds have been obtained using inflation.
Inflation is a procedure which works in any symplectic $4$-manifold $(M, \omega)$ and which consists
of altering the cohomology class of $\omega$ by adding to it  multiples of
the Thom class of a symplectic surface in $M$ of nonnegative self-intersection
(see Lalonde-McDuff \cite{Lalonde_McDuff1} and McDuff \cite{McDuff96, McDuff01}).

Inflation first appeared in the papers on the classification of symplectic structures
on ruled $4$-manifolds (see \cite{Lalonde_McDuff1, Lalonde_McDuff2}) where the authors use it in
a combination with another idea. We give a very rough sketch of this.
Let $(\bar X, \omega)$ be a ruled symplectic $4$-manifold, let $J\in \J$ be
an almost complex structure, let $\F$ be the
corresponding $J$-holomorphic foliation of $\bar X$ by spheres in the homology class $F$,
and let $\Sigma$ be a symplectic section of $\F$.
Consider other symplectic forms on $\bar X$ which have the form $\omega + k \pi^*\sigma$,
where $\pi: \bar X\rightarrow \Sigma$ denotes the projection along $\F$, $\sigma$ is an area form on $\Sigma$,
and $k \ge 0$ is any nonnegative constant.
It is easy to check that the forms $\omega_k = \omega + k \pi^*\sigma$ are indeed symplectic for any $k \ge 0$:
$\omega_k\wedge \omega_k = \omega^2 + 2k\omega\wedge \pi^*\sigma$
and both terms evaluate positively on a positively oriented basis of tangent vectors.
However the cohomology classes of $\omega_k$ are now different from the cohomology class of $\omega$.
Under certain conditions one can inflate along a suitable symplectic surface, obtaining
symplectic forms $\omega_k'$ whose cohomology classes are multiples of $[\omega]$.
Rescaling the symplectic forms, Moser's argument can be applied to produce interesting results.

In section $7$ we prove a version of the above argument which has
the advantages that no assumption on the self-intersection of
surfaces is made and that it applies nicely to families. This will
allow to prove the main theorem stated in the introduction. This
proof and its simple corollaries are discussed in section $8$.

\section{Holomorphic foliations}

Let $\J$ be the set of $\omega$-tame almost complex structures on $\bar X$.
For a symplectic surface $\Sigma\in \S_d(\bar X)$ we denote by
$\J_\Sigma \subset \J$ the subset of almost complex structures for which
$\Sigma$ is holomorphic.
For $J\in \J$ and a homology class $A\in H_2(M)$ we define
$\M(J; A)$ as the moduli space of irreducible $J$-holomorphic spheres representing
the homology class $A$. We abbreviate $\M(J; F)$ as $\M(J)$.

By \cite{HLS}, for any $J\in \J$,
$\M(J)$ is a smooth manifold. If nonempty, then by the index formula its dimension is $2$.
Since $F\cdot F = 0$, different spheres in $\M(J)$ cannot intersect by the
positivity of intersection.


\begin{proposition}\label{foliations_existence}
Let $(\bar X, \omega)$ be a ruled symplectic $4$-manifold. Let $J\in \J_\Sigma$.
\begin{itemize}
\item
Suppose that either $g\ge 1$, or that $g=0$, $d\ne 0$ and the
cohomology assumption $(*_d)$ is satisfied.
Then the spheres in $\M(J)$ form a holomorphic foliation of $\bar X$.
Moreover, each such sphere intersects the surface $\Sigma$ uniquely and
transversely. Furthermore, these foliations depend smoothly on $J$.
\item If $g=d=0$, then
the spheres in $\M(J; A_0)$ form a holomorphic foliation of $\bar X = S^2\times S^2$.
Moreover, each such foliation contains the sphere $\Sigma$ and these foliations depend smoothly on $J$.
\end{itemize}
\end{proposition}

\begin{proof}
Consider the first part of the proposition.
Given that the spheres in $\M(J)$ form a foliation of $\bar X$,
the second statement follows from the positivity of
intersection with $\Sigma$ and the third from the general theory of holomorphic curves
(see \cite{Gromov, McDuff_Salamon2}).
Now, the fact that for $J\in \J$ the $J$-holomorphic spheres in $\M(J)$ produce a
foliation of $\bar X$ is well-known when $g>0$,
or when $g=0$, $d\ne 0$ and $J\in \J$ is generic
(that is away from a subset of codimension $2$), see \cite{McDuff_Salamon2}.
Thus it remains to prove the first statement only
in the case $g=0$, $d\ne 0$ under the cohomology assumption $a_d > 0$ and $a_{-d} > 0$.

Assuming $g=0$, $d\ne 0$,
let $J_n\in \J$ be a sequence of generic almost complex structures approximating $J\in \J_\Sigma.$
By compactness theorem there is a (possibly) cusp $J$-holomorphic curve through every point of $\bar X$
which represents the homology class $F$.
We prove by contradiction that such a curve must in fact be nondegenerate.
We break the proof into two cases.

\vskip5pt
\noindent
{\bf Case $d>0$.}

\noindent
Suppose that $C = \cup_{i=1}^{n} C_i$ is a cusp $J$-holomorphic curve representing
the homology class $F$ and $n\ge 2$.
Write $[C_i] = k_i A_{-d} + m_i F$.
Thus $F = \sum_{i=1}^n (k_i A_{-d} + m_i F),$ with $\sum k_i = 0$ and $\sum m_i = 1$.
By the cohomology assumption, $\omega([\Sigma]) = a_d = a_{-d} + d\cdot f > f$ and so
no component $C_i$ can be a multiple cover of $\Sigma$. Thus by the positivity
of intersection with $\Sigma$, $m_i\ge 0$ for all $i$.
It follows that exactly one of the $m_i$'s is $1$ (we denote this $m_i$ by $m_1$)
and the rest are zero.
Also note that $k_1 < 0$, for otherwise $\omega([C_i]) > f$.
Thus the splitting in homology becomes
$F = (F - k A_{-d}) + \sum_{i=2}^n k_i A_{-d},$ where $k = -k_1 = \sum_{i=2}^n k_i > 0$
and $k_i \ge 0$ for $i\ge 2$.

We note that $F\cdot (F - kA_{-d}) = -k < 0$ and so there cannot
be any non-cusp $J$-holomorphic spheres representing the class $F$.
Choose a point $p$ not on $C$.
There must exist a cusp $J$-holomorphic sphere $C'$ through $p$ representing $F$.
Repeating previous arguments, $C'$ must be also split as $\cup_{i=1}^{n'} C_i'$,
with $[C_1'] = F -k'A_{-d}$, $k' > 0$.
However, $(F - kA_{-d}) \cdot (F - k'A_{-d}) = -k -k' - k k' d < 0$,
which is a contradiction since both $C_1$ and
$C_1'$ are $J$-holomorphic.

\vskip5pt
\noindent
{\bf Case $d<0$.}

\noindent
Again, suppose that $C = \cup C_i$ is a cusp
$J$-holomorphic curve in the homology class $F$.
Since $\Sigma^2 = d <0$, any holomorphic curve whose homology class is a multiple
of $A_d$ must be geometrically a multiple cover of $\Sigma$.
We group together those components which represent multiple covers of $\Sigma$ and those which do not.
Thus in homology
$F = \sum_{j=1}^n s_j A_d + \sum_{i=1}^{\tilde n}(k_i A_{-d} + m_i F )$
with $\tilde n > 0$ and no component in the second sum is a multiple cover of $\Sigma$.
By the positivity of intersections with $\Sigma$, $m_i \ge 0$ for all $i$.
By the cohomology assumption
$a_{-d} = a_d + |d| \cdot f$, in particular
$a_{-d} > f$.
Hence for every component in the second sum,
$m_i \ge 1$ and $k_i < 0$.
Clearly $s_j>0$ for all $j$.
Let $\bar k_i = -k_i > 0$.
Since $A_d = A_{-d} + dF$,
$\sum s_j = \sum \bar k_i$ and
$1 + \sum s_j d = \sum m_i$.

Since the area of each component must be positive,
$\omega(m_i F - \bar k_i A_{-d}) = m_i f - \bar k_i a_{-d} > 0.$
Recall that $a_{-d}/f > |d|$ and so $m_i \ge \bar k_i |d| +1 $.
Adding this inequalities for $i=1, \dots \tilde n$, we get
$$\sum {m_i} \ge \tilde n + \sum {\bar k_i |d| } = \tilde n + \sum s_j |d| = \tilde n -1 + \sum m_i.$$
It follows that $\tilde n - 1 \le 0$ and so the second sum consists of an exactly one term.
Let $k = \bar k_1$, $m= m_1$, then $\sum s_j = k > 0$ and $m_1 = 1 - kd$.
Thus in homology the splitting of $C$ must
be of the form
$F  = ( k \Sigma)  + ( (1 - kd ) F - k A_{-d} )$
with $(1 - kd ) F - k A_{-d}$ represented by a $J$-holomorphic sphere.

Note that $( (1 - kd ) F - k A_{-d}) \cdot F = -k < 0$, and so there cannot
be any non-cusp $J$-holomorphic spheres in the class $F$.
Choose a point $p$ not on $C$ and not on $\Sigma$,
and let $C'$ be a cusp $J$-holomorphic curve
through $p$ representing $F$.
By the above arguments $C'$ must contain a component $C_1'$ with homology
$ ( 1 - k' d) F - k' A_{-d}$ where $k' > 0$.
We note that
$C_1 \cdot C_1' = ((1 - kd ) F - k A_{-d} )\cdot ( ( 1 - k' d) F - k' A_{-d}) ) =
 -k ( 1 - k'd) - k'(1 - kd) - kk'd = -k - k' + kk'd < 0,$
contradicting to the positivity of intersection.

\vskip5pt
\noindent
Consider the second part of the proposition.
We note that since $\Sigma$ is $J$-holomorphic, there cannot be
any cusp $J$-holomorphic spheres representing the homology class $A_0$.
The existence of the required foliation
and its smooth dependence on $J$ now follow by standard methods.
\end{proof}

\section{Diffeomorphisms and smooth isotopy}

Assume that $(\bar X, \omega)$ is a ruled symplectic
$4$-manifold, that $\Sigma\subset \S_d(\bar X)$ is a symplectic surface, and that
the conclusions of the first part of proposition~\ref{foliations_existence} hold.
In other words, assume that for every $J\in \J_\Sigma$
the spheres in $\M(J)$ form a holomorphic foliation of $\bar X$.
This is automatically satisfied if one imposes the cohomology assumption $(*_d)$,
but this is also satisfied when $g=d=0$ and $f(\omega)\le a_0(\omega)$.

We will write $\F(J)$ for the
foliation given by $\M(J)$.
Note that we have a canonical diffeomorphism $\M(J)\rightarrow \Sigma$ associating to a holomorphic sphere
in $\M(J)$ its point of intersection with $\Sigma$.

Suppose that $J_t$, $t\in [0,1]$, is a family of almost complex
structures in $\J_\Sigma$. Let $\F_t = \F(J_t)$ be the family of
corresponding foliations. For $p\in \Sigma$ we denote by $C_{t,p}$
the sphere in $\F_t$ intersecting $\Sigma$ at $p$. We will define
a canonical vector field $Y_t$ generating a family $\Phi_t$ of
diffeomorphisms of $\bar X$ mapping the spheres $C_{0,p}$ to the
spheres $C_{t,p}$.
Denote by $N_{t,p}$ the symplectic normal bundle to $C_{t,p}$.
For each $t$ and $p$, for $|\epsilon|$ sufficiently small,
$J_{t+\epsilon}$-holomorphic spheres $C_{t,p}$
can be viewed as sections $\Gamma_{t,p,\epsilon}$ of $N_{t,p}$.
Define
$Y_t(p) = {\frac d {d\epsilon}}\bigl|_{\epsilon=0} \Gamma_{t,p,\epsilon}$
(note that $Y_t$ vanishes on $\Sigma$) and
let $\Phi_t$ be the time $t$ flow of $Y_t$, $0\le t\le 1$.
Thus $\Phi_t$ maps the spheres $C_{0,p}\in \F_0$ to the corresponding spheres $C_{t, p} \in \F_{t}$,
and $\Phi_t = id$ on $\Sigma$.

We are now ready to prove the smooth version of the main theorem (see theorem~\ref{intro_ruled}):
\begin{lemma}\label{smooth_contraction}
Let $(\bar X, \omega)$ be a ruled symplectic $4$-manifold.
Fix $d$ and suppose that the cohomology assumption $(*_d)$ holds.
Let $\Sigma\subset \S_d(\bar X, \omega)$ be a symplectic surface.
Let $\Sigma_\lambda \subset \S_{-d}(\bar X \setminus \Sigma, \omega)$
be a family of symplectic surfaces parameterized by $\lambda\in S^n$, $n\ge 0$.
Then the family $\Sigma_\lambda$ can be contracted through smooth surfaces
within $\bar X\setminus \Sigma$.
\end{lemma}

\begin{proof}
Given a family of symplectic surfaces $\Sigma_\lambda$, $\lambda\in S^n$, we can choose
a family $J_\lambda\in \J_\Sigma$ of almost complex structures such that
for each $\lambda$ the surface $\Sigma_\lambda$ is $J_\lambda$-holomorphic.

Since $\J_\Sigma$ is contractible, we can
contract the loop $J_\lambda$ to any given almost complex
structure $J_1\in \J_\Sigma$.
In other words, there exists a family $J_{\lambda,t}\subset \J_\Sigma$, $t\in [0,1]$, so that
\begin{enumerate}
\item
$J_{\lambda, 0} = J_\lambda$
\item
$J_{\lambda, 1} = J_1$.
\end{enumerate}
Since the assumptions of proposition~\ref{foliations_existence} are satisfied,
for each $\lambda$ and $t$ there exists a foliation $\F_{\lambda,t} = \F(J_{\lambda,t})$ of $\bar X$ by
holomorphic spheres in the homology class $F$. Moreover, $\F_{\lambda,t}$ depends smoothly on
the two parameters.
We let $\F_1 = \F(J_1)$ be the foliation corresponding to $J_1$.
Fixing $\lambda$, we define the canonical vector field $Y_{\lambda, t}$
whose time $t$ flow $\Phi_{\lambda, t}$ maps the foliation $\F_{\lambda, 0}$ to the
foliation $\F_{\lambda, t}$ and which is $id$ on $\Sigma$.
Note that by positivity of intersections, each surface $\Sigma_\lambda$ intersects every
sphere in $\F_{\lambda, 0}$ uniquely and transversely, and thus can be viewed as a
section of the foliation $\F_{\lambda, 0}$.
Define $\Sigma_{\lambda, t} = \Phi_{\lambda, t} (\Sigma_\lambda)$.
Then $\Sigma_{\lambda, t}$ a section of $\F_{\lambda, t}$ for every $\lambda$ and $t$, and
$\Sigma_{\lambda, t}$ provides a smooth isotopy between the $S^n$-family
$\Sigma_\lambda = \Sigma_{\lambda, 0}$ and an $S^n$-family of sections $\Sigma_{\lambda, 1}$ of $\F_1$.
By construction, all the surfaces $\Sigma_{\lambda, t}$ are disjoint from $\Sigma$.
Since the set of sections of $\F_1$ which are disjoint from $\Sigma$ is contractible,
the $S^n$-family $\Sigma_{\lambda, 1}$ can be contracted to any given section of $\F_1$
disjoint from $\Sigma$ through similar sections. This provides the required smooth contraction
of $\Sigma_\lambda$.

For the purposes of section $8$ we will still denote the extended isotopy from $\Sigma_\lambda$ to
a section of $\F_1$ by $\Sigma_{\lambda, t}$. We can extend the initial family $\F_{\lambda, t}$
of foliations by a constant family $\F_1$.
Since all the objects corresponding to $(\lambda, 1)$ are identified,
we can think of a pair $\mu = (\lambda, t)$ as of a point in $B^{n+1}$.
\end{proof}

\section{Straightening foliations}

We show how to modify symplectic foliations so that they intersect with
symplectic surfaces {\it nicely}.

\begin{lemma}\label{straighten_foliation}
Let $(\bar X, \omega)$ be a ruled symplectic $4$-manifold.
Let $\Sigma\subset \bar X$ be a symplectic surface
and fix a standard presentation (\ref{standard_form}) of $\omega$ near $\Sigma$.
Given any symplectic foliation of $\bar X$ intersecting $\Sigma$
positively and transversely, one
can modify it in a (arbitrarily small) neighborhood of $\Sigma$ so that the
new foliation coincides with the foliation by fibers of $\pi$ in some (even smaller) neighborhood
of $\Sigma$.
\end{lemma}

\begin{proof}
Fix a point $p\in \Sigma$. Denote by $C_p$ the curve in the foliation intersecting $\Sigma$ at $p$
and denote by by $F_p$ the fiber of $\pi$ through $p$ (which is defined in a neighborhood of $\Sigma$).

\vskip5pt
We first consider the case $d=\Sigma\cdot \Sigma = 0$.
In a neighborhood $U\times V$
of $p$ the symplectic form can be written as
$\omega = dx_1\wedge dx_2 + dy_1\wedge dy_2 = dx_1\wedge dx_2 + rdr\wedge d\theta$,
where $(x_1, x_2)\in U\subset \Sigma$ are coordinates on a neighborhood $U$ of $p$ on $\Sigma$,
$(y_1, y_2)\in V\subset F_p$ are rectilinear coordinates on a neighborhood $V$ of $p$ on $F_p$, and $(r, \theta)$ are
the corresponding polar coordinates on $V$.

By restricting to a possibly smaller neighborhood $U\times V$ we can assume that $C_p \cap (U\times V)$
is a graph over $V$, thus in $(r, \theta, x_1, x_2)$-coordinates $C_p$ has the form
$C_p = \{ (r, \theta, x_1(r, \theta), x_2(r, \theta)) \  | \ (r,\theta)\in V\}$.
Let $q\in C_p\cap (U\times V)$. Notice that since $C_p$ intersects $\Sigma$ positively and
transversely, the restriction of the form $rdr\wedge d\theta$ to $T_qC_p$ induces positive
orientation.
The tangent space $T_qC_p$ (when $q\ne p$) is spanned by
$e_1 = (1, 0, \frac {\partial x_1} {\partial r}(r, \theta), \frac {\partial x_2} {\partial r}(r, \theta) )$
and $e_2 = (0, 1, \frac {\partial x_1} {\partial \theta}(r, \theta), \frac {\partial x_2} {\partial \theta}(r, \theta))$.
Since $T_qC_p$ is symplectic, for $r>0$ we have that
\begin{equation}
\omega(e_1, e_2) = r + D(r, \theta) > 0 ,
\end{equation}
where $D(r, \theta) =
[\frac {\partial x_1} {\partial r}\frac {\partial x_2} {\partial \theta} -
\frac {\partial x_2} {\partial r}\frac {\partial x_1} {\partial \theta}](r, \theta)$.
Since $C_p$ is smooth and symplectic at $p$, we have that $\lim_{(r,\theta)\rightarrow 0} ( 1 + D(r, \theta)/r ) $ exists
and is positive.
Since $V$ is compact, there is a constant $\delta > 0$ so that
\begin{equation}
(1- \delta) r +  D(r, \theta) > 0
\end{equation}
for $r>0$.

Given $\epsilon_1>0$ sufficiently small,
and given $\epsilon_0>0$ sufficiently small compared to $\epsilon_1$, one can choose a nonnegative
nondecreasing function $\beta(r)$ so that $\beta(r) = 0$ for $r\le \epsilon_0$,
$\beta (r) =  r$ for $ r\ge \epsilon_1$, $\beta(r)\le r$ and $(1-\delta )\beta'( r) < 1$.

We choose $\epsilon_1>0$ sufficiently small so that $\{r\le \epsilon_1\} \subset V$ and so that it
satisfies whatever smallness assumptions of the lemma. We choose $\epsilon_0$ and $\beta(r)$ as in
the preceding paragraph.
We define the curve $\tilde C_{p}$ by ``radially-rescaling" $C_p$:
$\tilde C_p = \{ (r, \theta, x_1(\beta(r),  \theta), x_2(\beta (r), \theta)) \ |\ (r, \theta) \in V \}$.
We note that $\tilde C_{p}$ defined in this way is a smooth sphere which coincides
with $F_{p}$ for $ r\le \epsilon_0$ and which coincides with $C_{p}$ for $ r\ge \epsilon_1$.
We claim that  $\tilde C_p$ is $\omega$-symplectic. We only need to verify this for
$\epsilon_0 < r < \epsilon_1$. For $q\in \tilde C_p$ lying above $(r, \theta)\in F_p$,
the tangent plane $T_q\tilde C_p$ is spanned by
$\tilde e_1 = (1, 0, \beta'(r) \frac {\partial x_1} {\partial r}(\beta(r), \theta), \beta'(r) \frac {\partial x_2} {\partial r}(\beta(r), \theta) )$
and
$\tilde e_2 = (0, 1, \frac {\partial x_1} {\partial \theta}(\beta(r), \theta), \frac {\partial x_2} {\partial \theta}(\beta(r), \theta))$.
Thus at $q$ we have:
$$\omega_q(\tilde e_1, \tilde e_2) = r + D(\beta(r), \theta) \beta'(r) \ge r - (1-\delta) r \beta'(r) > 0,$$
verifying that $\tilde C_p$ is symplectic (at $q$).

Next, for each $p'\in \Sigma$ we let $C_{p'}$ be the curve in the foliation through $p'$ and
we let $F_{p'}$ be the fiber through $p'$.
We note that for $p'$ sufficiently close to $p$, $C_{p'}$ is a graph over $F_{p'}\cap (U\times V)$.
Moreover (being symplectic is an open condition), for $p'$ sufficiently close to $p$
the radially-rescaled curves $\tilde C_{p'}$ (defined in the analogous way with the same function $\beta(r)$)
are also $\omega$-symplectic.
By a compactness argument one can choose the same $\beta$ to work for every sphere in the foliation.
Thus the ``radially-rescaled" spheres $\tilde C_{p}$ depend
smoothly on $p$ and satisfy all the requirements of the lemma.

\vskip5pt
The proof is almost analogous in the case $d\ne 0$.
Since $T_p\Sigma$ and $T_pF_p$ are $\omega$-orthogonal,
$\omega$ has the form $\omega_p = dx_1\wedge dx_2 + rdr\wedge d\theta$
at the point $p$ and is sufficiently close to this form in a sufficiently small neighborhood $U\times V$ of $p$.
More precisely, for $q\in C_p$ in a sufficiently small neighborhood $U\times V$ of $p$ the planes
$T_qC_p$ are $\omega_p$-symplectic. By the previous construction we can define an $\omega_p$-symplectic
``radially-rescaled" curve $\tilde C_p$.
The tangent planes to $\tilde C_p$ are $\omega_q$-symplectic for all $q$ sufficiently close to $p$.
Thus shrinking the defining function $\beta(r)$ if necessary
(i.e. considering $\frac 1 K \beta(Kr)$ for $K$ sufficiently large), we can obtain the required $\omega$-symplectic
rescaling of $C_p$.
The rest of the proof follows by a compactness argument.
\end{proof}

\begin{remark}
Suppose that $\Sigma_\lambda$, $\lambda\in S^n$, is a family of $\omega$-symplectic surfaces in $\bar X$,
and $\F_\lambda$ is a family of $\omega$-symplectic foliations of $\bar X$ such that for each $\lambda$ every
sphere in $\F_\lambda$ intersects $\Sigma_\lambda$ positively and transversely.
We fix a family of standard presentations of $\omega$ near $\Sigma_\lambda$ and denote by
$\pi_\lambda$ the corresponding family of projection maps (each defined in a neighborhood of the
corresponding surface).
Clearly one can choose the same function $\beta(r)$ for every sphere in every foliation.
Thus we can obtain a new family $\F_\lambda$ of $\omega$-symplectic foliations
so that for each $\lambda$ the foliation $\F_\lambda$ intersects $\Sigma_\lambda$ nicely.
\end{remark}

\section {An inflation-like argument}

The heart of this paper lies in the following simple argument.

\begin{proposition}\label{inflation}
Let $(\bar X, \omega)$ be a ruled symplectic $4$-manifold. Let
$\Sigma\subset \S_d(\bar X)$ be a symplectic surface and fix a
standard presentation (\ref{standard_form}) of $\omega$ near
$\Sigma$. By scaling $\omega$ we may assume that $a_d(\omega) = |d|$
if $d\ne 0$ and $a_d(\omega) = 1$ if $d=0$.
If $d>0$ also assume that the cohomology assumption $(*_d)$ is satisfied.
Let $\F$ be a symplectic foliation of $\bar X$ nicely
intersecting $\Sigma$. Then there exists a family $\hat \omega_k$,
$k\ge 0$, of cohomologous symplectic forms on $\bar X$ which
smoothly depend on the parameter $k$ and which satisfy:
\begin{itemize}
\item $\hat \omega_k = \omega$ near $\Sigma$,
\item $\hat \omega_k$ is a multiple of $\omega + k\pi^*\sigma$ outside $U$,
\item$\hat \omega_k = \omega$ when $k=0$.
\end{itemize}
\end{proposition}

\begin{proof}

Recall that for $k \ge 0$ the forms
$\omega_k \Def \omega + k \pi^*\sigma$ are symplectic.
We break the construction into $3$ cases.

\vskip5pt
\noindent
{\bf Case: $d>0$:}

\noindent
Define the form $\omega_k'$ by rescaling the form $\omega_k$:
$$ \omega_k' =\lambda_k \omega_k,$$
where $$ \lambda_k = \frac {1 - f }  {1+k - f }. $$
Note that $\lambda_k\in (0,1]$ and $a_{-d}(\omega_k') = a_{-d}(\omega)$.

Choose $\epsilon$ so that $\{r\le \epsilon\} \subset U$.
Let $\beta(r)$ be a positive decreasing function which equals $1-r^2$ near $r=0$
and equals $\frac {1-f} {k+1-f} (k + 1 - r^2)$ near $\epsilon$.
To verify that such a function exists, notice that when $k=0$ the expressions
$\frac {1-f} {k+1-f} (k + 1 -r^2)$ and $1-r^2$ coincide, and when $k>0$ the
expression $1-r^2$ is strictly greater than $\frac {1-f} {k+1-f} (k + 1 -r^2)$,
at least for all $r$ sufficiently small.
The function $\beta(r)$ can be constructed by an appropriate smoothing between these
two functions.

Define the form $\hat \omega_k = d(\beta(r)\alpha )$. By construction it is a symplectic form,
it equals $\omega$ near $r=0$ and equals $\lambda_k (\omega + \pi^*\sigma )$ outside $U$.
Since $H_2(\bar X; \R)$ is generated by $A_d$ and $A_{-d}$
and since $a_d(\hat \omega_k) = a_d(\omega)$ and $a_{-d}(\hat \omega_k) = a_{-d}(\omega)$,
the form $\hat \omega_k$ is cohomologous to $\omega$. Finally note that the construction can be done
to depend smoothly on $k$.

\vskip5pt
\noindent
{\bf Case: $d<0$:}

\noindent
Again define
$ \omega_k' =\lambda_k \omega_k$
with
$$ \lambda_k = \frac {1 + f }  { 1+k + f }.$$
Note that $\lambda_k\in (0,1]$ and $a_{-d}(\omega_k') = a_{-d}(\omega)$.

Choose $\epsilon$ so that $\{r\le \epsilon\} \subset U$.
Let $\beta(r)$ be a positive increasing function which equals $1+r^2$ near $r=0$
and equals $\frac {1+f} {k+1+f} (k + 1 + r^2)$ near $\epsilon$.
To verify that such a function exists, notice that when $k=0$ the expressions
$\frac {1+f} {k+1+f} (k + 1 + r^2)$ and $1+r^2$ coincide, and when $k>0$ the
expression $1+r^2$ is strictly smaller than $\frac {1+f} {k+1+f} (k + 1 + r^2)$,
at least for all $r$ sufficiently small.
The function $\beta(r)$ can be constructed by an appropriate smoothing between these
two functions.

Define the form $\hat \omega_k = d(\beta(r)\alpha )$. By construction it is a symplectic form,
it equals $\omega$ near $r=0$ and equals $\lambda_k (\omega + \pi^*\sigma )$ outside $U$.
Since $a_d(\hat \omega_k) = a_d(\omega)$ and $a_{-d}(\hat \omega_k) = a_{-d}(\omega)$,
the form $\hat \omega_k$ is cohomologous to $\omega$.
Again the construction can be done to depend smoothly on $k$.

\vskip5pt
\noindent
{\bf Case: $d=0$:}

\noindent
With $\lambda_k = \frac 1 {k+1}$ define
$$\omega_k' = \lambda_k \omega_k = \pi^*\sigma + \lambda_k dy_1\wedge dy_2.$$
Note that $\lambda_k\in (0,1]$ and $a_0(\omega_k') = a_0(\omega)$.

Choose $\epsilon$ so that $\{r\le \epsilon\} \subset U$.
Choose a positive function $\phi(r)$ so that
$\phi(r) = 1$ near $r=0$, $\phi(r) = \lambda_k$ for $r\ge \epsilon$ and
$f( \pi^*\sigma + \phi(r)dy_1\wedge dy_2  ) = f(\omega)$.
Define $\hat \omega_k  = \pi^*\sigma + \phi(r)dy_1\wedge dy_2 $.
This is a symplectic form, it equals $\omega$ near $r=0$, equals $\lambda_k (\omega + k\pi^*\sigma)$
outside $U$, and is cohomologous to $\omega$.
The construction can be done to depend smoothly on $k$.

\end{proof}

\section{Topology of spaces of surfaces}

We can now prove the main theorem (see theorem~\ref{intro_ruled})
if either $d\ne 0$ or $g\ne 0$.

\begin{theorem}\label{main_theorem_ruled}
Let $(\bar X, \omega)$ be a ruled symplectic $4$-manifold.
Fix $d\in \Z$ and suppose that the cohomology assumption $(*_d)$ holds.
Let $Z\in \S_d(\bar X)$ be a symplectic surface representing the homology class $A_d$.
Then the set $\S_{-d}(\bar X\setminus Z)$
of symplectic surfaces representing the homology class $A_{-d}$
and disjoint from $Z$ is non-empty and contractible.
\end{theorem}

\begin{proof}

First we show that $\S_{-d}(\bar X\setminus Z)$ is non-empty,
i.e. that we can find a symplectic surface $\Sigma\in \S_{-d}$
which is disjoint from $Z$.
In the following we will assume that
$d\ne 0$ (if $d=0$
then the statement is clearly true by moving $Z$ slightly off itself).

Choose an $\omega$-tame almost complex structure $J$ for which $Z$ is holomorphic.
By proposition~\ref{foliations_existence} we obtain
a $J$-holomorphic foliation $\F$ of $\bar X$ by the spheres in the homology class $F$,
each of which intersects $Z$ uniquely, transversely and positively.
We fix a standard presentation (\ref{standard_form}) of $\omega$ near $Z$.
By lemma~\ref{straighten_foliation} we can modify the foliation $\F$ in a neighborhood of $Z$
so that it remains symplectic and intersects $Z$ nicely (that is,
near $\Sigma$ coincides with the foliation given by
the fibers of $\pi$).

Choose any smooth section $S$ of $\F$ which is disjoint from $Z$. For $K\ge 0$
sufficiently large, $S$ is $\omega + K\pi^*\sigma$-symplectic. By proposition~\ref{inflation} there exists
a family $\hat \omega_k$ of symplectic forms which are cohomologous to $\omega$, coincide with $\omega$ near $Z$,
and so that $\hat \omega_0 = \omega$ and $S$ is $\hat \omega_K$-symplectic.

Consider the family $\hat \omega_k$ of symplectic forms, $k\in [0 ,K]$. By Moser's method there is a
symplectomorphism $\psi: (\bar X, \hat \omega_0 = \omega) \rightarrow (\bar X, \hat \omega_K)$
which is $id$ near $Z$. It follows that $\Sigma \Def \psi^{-1}(S)$ is an $\omega$-symplectic
surface disjoint from $Z$.

\vskip5pt

Next we show that $\S_{-d}(\bar X\setminus Z)$ is contractible.
Consider any $S^n$-family of symplectic surfaces $\Sigma_\lambda\in \bar X\setminus Z$.
In the proof of lemma~\ref{smooth_contraction}
we have constructed a family $\F_\mu$ of
symplectic foliations and a corresponding family $\Sigma_\mu$ of sections, so that:
$\Sigma_\mu = \Sigma_\lambda$ is $\omega$-symplectic for $\mu \in S^n = \partial B^{n+1}$,
$\Sigma_\mu$ is disjoint from $Z$ for all $\mu$, and
$\F_\mu$ intersects $Z$ uniquely, transversely and positively for all $\mu$.
Our goal is to modify the smooth contraction $\Sigma_\mu$
in $\bar X\setminus Z$ to a symplectic one.

We can assume that none of the surfaces $\Sigma_\mu$ intersect a sufficiently small neighborhood
of $Z$. Fixing a standard presentation (\ref{standard_form}) of $\omega$ near $Z$,
by the remark following lemma~\ref{straighten_foliation}
we can modify all the foliations $\F_\mu$ near $Z$ so that they intersect $Z$ nicely.
We denote by $\pi_\mu: \bar X\rightarrow Z$ the projection along the fibers of $\F_\mu$,
that is $\pi_{\mu}$ sends each sphere in the
foliation $\F_\mu$ to its intersection with $Z$.

For each $\mu\in B^{n+1}$ we can
choose $K_\mu$ sufficiently large so that $\Sigma_\mu$ becomes symplectic
with respect to the symplectic form $\omega_\mu' \Def \omega + K_\mu \pi_\mu^*\sigma$.
The values for $K_\mu$ can be chosen to depend smoothly on $\mu$ and
we can assume that $K_\mu = 0$ for $\mu\in \partial B^{n+1}$
(so that $\omega_\mu' = \omega$ for $\mu\in \partial B^{n+1}$).

Applying proposition~\ref{inflation}, for each $\mu$ we can find
a family $\hat \omega_{\mu, k}$ of symplectic forms which are cohomologous to $\omega$,
coincide with $\omega$ near $Z$, so that $\hat \omega_{\mu, 0} = \omega$ and
$\Sigma_\mu$ is $\hat \omega_{\mu, K_\mu}$-symplectic.

For each $\mu\in B^{n+1}$ we can apply Moser's argument to the family $\hat \omega_{\mu, k}$ for $k\in [0, K_\mu]$
to obtain a symplectomorphism $\psi_\mu : (\bar X, \omega) \rightarrow (\bar X, \hat \omega_{\mu, K_\mu})$.
Note that each $\psi_\mu$ is $id$ near $Z$ and $\psi_\mu = id$ when $\mu \in \partial B^{n+1}$.

The preceding construction can be made to depend smoothly on $\mu$.
It follows that $\tilde \Sigma_\mu = \psi_{\mu}^{-1}(\Sigma_\mu)$ is
a family of $\omega$-symplectic surfaces in $\bar X$
providing a symplectic contraction of $\Sigma_\lambda = \Sigma_\mu$, $\mu\in\partial B^{n+1}$,
in $\bar X\setminus Z$.

\end{proof}

Next consider the only case of theorem~\ref{intro_ruled} not covered by the theorem above, i.e.
the case $g=d=0$.

\begin{proposition}
Let $Z\in \S_0(S^2\times S^2, \omega)$ be a symplectic surface representing the
homology class $A_0$.
Then the set $\S_0(S^2\times S^2 \setminus Z, \omega)$ of
symplectic surfaces representing the homology class $A_0$ and disjoint from $Z$
is non-empty and contractible.
\end{proposition}

\begin{proof}
First note that $Z\subset S^2\times S^2$ has a standard symplectic neighborhood of the form
$Z\times D^2$ with a split symplectic form. Taking a constant section $Z\times \{pt\}$
we see that $\S_0(S^2\times S^2 \setminus Z, \omega)$ is non-empty.

Next, let $\lambda \in S^n$ parameterize a family in $S^2\times S^2 \setminus Z$.
 In other words, let
$\Sigma_\lambda$ be an $S^n$-family of embedded $\omega$-symplectic surfaces
homologous to $Z$ and disjoint from $Z$. One can choose
a family $J_\lambda\subset \J_Z$ of almost complex structures on $S^2\times S^2$
such that $\Sigma_\lambda$ is $J_\lambda$-holomorphic for each $\lambda$.
Since $\J_Z$ is
contractible, we can extend $J_\lambda$ to a $B^{n+1}$-family of almost complex
structures in $\J_Z$ which coincides with the original family for
$\lambda \in \partial B$.
According to proposition~\ref{foliations_existence}, for each $\lambda\in B^{n+1}$ there is a
$J_\lambda$-holomorphic foliation of $S^2\times S^2$ by spheres
representing the homology class $A_0$. We will denote
by $\B_{J_\lambda}$ the corresponding foliation of $S^2\times S^2 \setminus Z$
obtained by excluding $Z$.
Let's define
$$Y = \{ (\lambda, Z) \vert \lambda \in B^{n+1}, Z \in \B_{J_\lambda} \}.$$
This is a locally trivial fibration over $B^{n+1}$ (the projection map
is given by $(\lambda, Z) \rightarrow \lambda$) with the fiber homeomorphic to $D^2$.
Thus we can extend the section $\Sigma_\lambda$ defined for $\lambda\in \partial B^{n+1}$
to a global section of $Y$. This gives the required symplectic isotopy of $\Sigma_\lambda$
to a constant family.
\end{proof}

The local isotopy claimed in the introduction (see theorem~\ref{intro_local}) is now an easy consequence.
\begin{corollary}\label{main_theorem_local}
Let $(M, \omega)$ be a symplectic $4$-manifold and let $\Sigma\subset M$ be an embedded
symplectic surface. Then there is a (arbitrarily small) neighborhood $U$ of $\Sigma$
so that the set $\S(U, \omega)$ of embedded symplectic surfaces in $U$ homologous to $\Sigma$
is weakly contractible.
\end{corollary}
\begin{proof}
Let $U$ be any open standard symplectic neighborhood of $\Sigma$, i.e.
a multiple $c\omega$ of $\omega$ is
given by the standard form (\ref{standard_form}).
Let $\Sigma_\lambda$ be any $S^n$-family in  $\S(U, \omega)$.
As described in section $2$
we can compactify $(U, c\omega)$ to a standard symplectic sphere bundle $\bar X$.
In particular $U$ is symplectomorphic to $\bar X\setminus \Sigma_\infty$ for a
symplectic surface $\Sigma_\infty \in \bar X$. By theorem~\ref{main_theorem_ruled}
$\Sigma_\lambda$ can be symplectically contracted in $\bar X\setminus \Sigma_\infty$ and hence in $U$.
\end{proof}

The proof of the main theorem~\ref{main_theorem_ruled} can be easily generalized as follows:

\begin{proposition}
Let $(\bar X, \omega)$ be a ruled symplectic $4$-manifold.
Fix $d\ne 0$ and suppose that the cohomology assumption $(*_d)$ holds.
Let $Z_\lambda\in \S_d(\bar X)$, $\lambda\in S^n$, be a family of symplectic surfaces
representing the homology class $A_d$.
Then there exists a corresponding family $\Sigma_\lambda$ of symplectic surfaces
representing the homology class $A_{-d}$ so that for each $\lambda\in S^n$
the surface $\Sigma_\lambda$ is disjoint from $Z_\lambda$.
Moreover, the set of such families $\{\Sigma_\lambda\}$ is contractible.
\end{proposition}

\begin{corollary}
Let $(\bar X, \omega)$ be a ruled symplectic $4$-manifold satisfying the cohomology assumption $(*_d)$.
Then $\S_d(\bar X, \omega)$ is homotopy equivalent to $\S_{-d}(\bar X, \omega)$.
\end{corollary}

\begin{proof}
Consider the set $Y$ of pairs $(Z, \Sigma)$ with $Z\in \S_d$, $\Sigma\in \S_{-d}$,
and $Z$ and $\Sigma$ disjoint.
The projection maps $Y\rightarrow \S_d$ and $Y\rightarrow \S_{-d}$ are fibrations with contractible fiber
and hence are both homotopy equivalences.
\end{proof}

\section{Relations with groups of symplectomorphisms}
Motivated by \cite{Coffey}, we consider some obvious connections between the groups of symplectomorphisms
of a ruled symplectic $4$-manifold and the sets of symplectic surfaces.
Let $(\bar X, \omega)$ be a ruled symplectic $4$-manifold. Fix $d\ne 0$ and suppose that
$(\bar X, \omega)$ satisfies the homology assumption $(*_d)$.

First we introduce some notations.
We denote by $\Symp^h(\bar X, \omega)$ the group of symplectomorphisms of $\bar X$ which act
as identity on homology. Given $Z\in \S_d(\bar X)$, we denote by
$\Symp^h_{[Z]}(\bar X, \omega) \subset \Symp^h_{\bar Z}(\bar X, \omega)
\subset \Symp^h_Z(\bar X, \omega) \subset \Symp^h(\bar X, \omega)$
the subsets of $\Symp^h(\bar X, \omega)$ consisting of symplectomorphisms
which, respectively, fix a neighborhood of $Z$, fix $Z$, and preserve $Z$.

\vskip5pt

Consider the action of $\Symp^h_{[Z]}(\bar X, \omega)$ on $S_{-d}(\bar X\setminus Z)$.
By theorem~\ref{main_theorem_ruled}, any two symplectic surfaces in
$S_{-d}(\bar X\setminus Z)$ can be joined by a path of such symplectic surfaces.
None of these surfaces intersect a sufficiently small neighborhood of $Z$ and so
this path can be generated by a Hamiltonian which is identity in a neighborhood of $Z$.
It follows that the action above is transitive.
By the same theorem, the set $S_{-d}(\bar X\setminus Z)$ is in fact contractible.
The stabilizer of a surface $\Sigma\in S_{-d}(\bar X\setminus Z)$
is the set of symplectomorphisms of $\bar X$ which preserve $\Sigma$ and
fix a neighborhood of $Z$. As shown in \cite{Coffey}, this set is also contractible.
It follows that
\begin{lemma}
Assuming the cohomology assumption $(*_d)$,
the group of symplectomorphisms of $(\bar X, \omega)$ which act as identity on
homology and fix a neighborhood of a symplectic surface $Z\in \S_d(\bar X)$ is
contractible.
\end{lemma}

Next consider the action of $\Symp^h_Z(\bar X, \omega)$ on $S_{-d}(\bar X\setminus Z)$.
As before this action is transitive and the set $S_{-d}(\bar X\setminus Z)$ is contractible.
The stabilizer of a surface $\Sigma\in S_{-d}(\bar X\setminus Z)$ consists of
symplectomorphisms which preserve both $\Sigma$ and $Z$. By \cite{Coffey}, this set is
homotopy equivalent
to the set $D_{0, \infty}$ of fiber preserving diffeomorphisms
of a standard sphere bundle of degree $d$ over $\Sigma$ which also preserve the sections
$\Sigma_0$ and $\Sigma_\infty$.
It follows that $\Symp^h_Z(\bar X, \omega)$ is homotopy equivalent to $D_{0, \infty}$.

\vskip5pt

By \cite{Coffey}, given any two pairs of non-intersecting symplectic sections $(Z_1, \bar Z_1)$,
$(Z_2, \bar Z_2)$ with $Z_1, Z_2 \in \S_d(\bar X)$ and  $\bar Z_1, \bar Z_2 \in \S_{-d}(\bar X)$,
there exists a symplectomorphism in
 $\Symp^h(\bar X, \omega)$ mapping $Z_1$ to $Z_2$ and mapping
$\bar Z_1$ to $\bar Z_2$. In particular, the group $\Symp^h(\bar X, \omega)$ acts transitively on $\S_d(\bar X)$.
The stabilizer of this action corresponding to $Z\in S_d(\bar X)$ is precisely the set
$\Symp^h_Z(\bar X, \omega)$ of symplectomorphisms which preserve $Z$ considered before.
Therefore we obtain the following:
\begin{lemma}
Assuming the cohomology assumption $(*_d)$,
the group 
$\Symp^h(\bar X, \omega)$
acts transitively on the set $\S_d(\bar X, \omega)$
with the stabilizer homotopy equivalent to
$D_{0, \infty}$.
\end{lemma}

It is interesting to note the connection with the theorem
of \cite{Abreu_McDuff} (see also \cite{McDuff01}) which can be formulated as follows.
Let $\bar X = \Sigma\times S^2$ and let $\omega_\mu = \mu \sigma_\Sigma + \sigma_{S^2}$.
Denote by $\Symp_{id}(\bar X, \omega_\mu) \subset \Symp^h(\bar X, \omega_\mu)$ the group
of symplectomorphisms which are isotopic to $id$ as diffeomorphisms, and denote by
$D$ the identity component of the group of fiber preserving diffeomorphisms of $\Sigma \times S^2$.
The theorem states that as $\mu$ tends to infinity, the groups $\Symp_{id}(\bar X, \omega_\mu)$
tend to a limit which is homotopy equivalent to $D$.

\bibliography{isotopy23}
\nocite*

\vspace{0.1in}

Richard Hind\\
Department of Mathematics\\
University of Notre Dame\\
Notre Dame, IN 46556\\
email: hind.1@nd.edu

\vspace{0.1in}

Alexander Ivrii\\
D\'epartment de Math\'ematiques et de Statistique\\
Universit\'e de Montr\'eal\\
CP 6128 succ Centre-Ville\\
Montr\'eal, QC H3C 3J7, Canada\\
email: ivrii@DMS.UMontreal.CA

\end{document}